\documentclass[a4paper, 12pt]{article}
\usepackage[a4paper]{geometry}
\usepackage{amsmath,amssymb,amsthm}
\usepackage[ruled, boxed, linesnumbered]{algorithm2e}
\usepackage{esvect}
\usepackage{xcolor}
\usepackage{enumerate}
\usepackage{caption}
\usepackage{natbib}
\usepackage{latexsym}
\usepackage{graphics}
\usepackage{graphicx}
\usepackage{wrapfig}
\usepackage{lineno}
\usepackage{wasysym}
\usepackage{stmaryrd}
\usepackage[english]{babel}
\usepackage{hyperref} 
\usepackage{relsize}
\usepackage{authblk}
 \allowdisplaybreaks
\author[1]{Martin Ritchie}
\author[2]{Luc Berthouze}
\author[1]{Istvan Z. Kiss \thanks{Corresponding author: I.Z.Kiss@ssusex.ac.uk}}
\affil[1]{School of Mathematical and Physical Sciences, Department of Mathematics, University of Sussex, Falmer, Brighton BN1 9QH, UK}
\affil[2]{Centre for Computational Neuroscience and Robotics, University of Sussex, Falmer,  Brighton, BN1 9QH, UK} 
            
\begin{document}

\title{\bf Beyond clustering: Mean-field dynamics on networks with arbitrary
subgraph composition \\}

\maketitle

\begin{abstract}
Clustering is the propensity of nodes that share a common neighbour to be connected. It is ubiquitous in many networks but poses many modelling challenges. Clustering typically manifests itself by a higher than expected frequency of triangles, and this has led to the principle of constructing networks from such building blocks. This approach has been generalised to networks being constructed from a set of more exotic subgraphs. As long as these are fully connected, it is then possible to derive mean-field models that approximate epidemic dynamics well. However, there are virtually no results for non-fully connected subgraphs. In this paper, we provide a general and automated approach to deriving a set of ordinary differential equations, or mean-field model, that describes, to a high degree of accuracy, the expected values of system-level quantities, such as the prevalence of infection. Our approach offers a previously unattainable degree of control over the arrangement of subgraphs and network characteristics such as classical node degree, variance and clustering. The combination of these features makes it possible to generate families of networks with different subgraph compositions while keeping classical network metrics constant. 
Using our approach, we show that higher-order structure realised either through the introduction of loops of different sizes or by generating clustered networks based on different subgraphs, leads to significant differences in epidemic dynamics despite controlling for basic network metrics.
\end{abstract}

\section{Introduction}\label{intro}

Network models have revolutionised our way of thinking about complex phenomena such as the spreading of disease, information transmission and processing in the brain, and the formation and interaction of social groups. Mathematical epidemiology, in particular, has embraced and benefited greatly from the use of networks as a modelling paradigm, with examples ranging from data-driven models \citep{tildesley2010impact,barabasi1999emergence,kiss2006network} to theoretical models \citep{newman2002spread,pastor2001epidemic,Keeling22041999,sharkey2013exact}. These have been used to study the impact of different network properties on how diseases break out and spread. Network models have led to greater clarity in understanding and quantifying the impact of contact heterogeneity, preferential mixing and community structure, including households \citep{ball2001stochastic,ball2010analysis}. Although clustering of contacts or transitivity, i.e., the propensity of nodes with a common neighbour to be connected, is pervasive in many real-world networks, it continues to pose many significant challenges to the community, both from the viewpoint of network generation and, even more so, from that of deriving well-performing approximate models.

To investigate the impact of network properties, one can either use empirical networks or synthetic ones that have been generated from theoretical network models with tunable properties \citep{newman2002spread,karrerclust2010,molloy1995critical}. Many algorithms exist for clustering, but it is generally the case that focusing on achieving a particular clustering leads to changes above and beyond those controlled by the algorithm. This can preclude correct analysis of the impact of clustering \citep{karrerclust2010,ritchie2014higher,house2009motif,housetracing,
Keeling22041999,green2010large,milo2002network,colomer2013deciphering,
miller2009percolation,gleeson2010clustering,kiss2008comment}. When looking at the impact of higher-order structure, for example, it is important that the degree distribution remain the same between networks with different clustering. Some algorithms in this direction have been proposed \citep{karrerclust2010,vmclust,miller2009percolation,newman2009random} and are based on the notion of subgraphs, where clustering is achieved by mixing fully-connected subgraph types, such as fully-connected triples or quadruples, and non-fully connected subgraphs, such as overlapping triangles. Using such networks, Volz et al. \citep{vmclust} have developed a low-dimensional ODE model that approximates well the expected value of a number of system-level quantities, and Karrer \& Newman \citep{karrerclust2010} have provided final epidemic size results for networks built by using different mixtures of subgraphs. Furthermore, House and colleagues \citep{house2009motif,wrap5586} generalised the pairwise approach to closure at the level of all possible subgraphs involving four nodes.  However, a number of outstanding issues remain. The Volz et al. model, which provides time evolution, can handle well only fully-connected subgraphs. Karrer \& Newman's approach, which combines a wider variety of subgraphs, can only characterise large-time limits. Finally, to our knowledge, House et al.'s approach has not been compared to stochastic simulations and it will perform poorly for heterogeneous networks \citep{house2009motif}.

In this paper, we provide a general and automated approach to deriving a set of ODEs that describe, to a high degree of accuracy, the expected values of prevalence or number of recovered individuals for networks that are generated based on an arbitrary set of subgraphs. This is achieved by a rigorous separation of the role of nodes within the subgraphs and by using the probability generating function (PGF) formalism to correctly track: (a) the distribution of subgraphs to which nodes belong and (b) the excess degree that is generalised from the classical notion of a stub of a single edge to different corner types given by subgraphs. This is a significant step forward as it allows us to: (a) accurately model and analyse dynamical processes on networks with higher-order structure, thus increasing model realism, (b) map out the impact of clustering in the classic sense, and more importantly, its impact at a higher level involving four or more nodes \citep{ritchie2014higher}, and (c) provide much needed insights into the role of small subgraphs or network motifs/units in epidemiology and systems biology.

The paper is organised as follows. We first review how the probability generating function (PGFs) can be used to derive ODEs that capture epidemic dynamics on configuration model (CM) networks. Such PGF-based models operate by using the versatile properties of the PGF whereby it allows us to keep track of the fraction of susceptible individuals, their degree and excess degree. Next, we generalise the CM to the \emph{hyperstub configuration model} (HCM). The HCM is a network construction algorithm that selects and connects hyperstubs as prescribed by the building blocks or subgraphs of the network, rather than at random. With a basic understanding of both the network and epidemic models, we then generalise the PGF formalism to HCM networks. This section includes a step-by-step explanation of the model derivation with  examples for a particular network and a detailed presentation of the code-generating algorithm. A key component of the generalised model is to label and track the position of each and every node in all subgraphs in order to avoid any ambiguity as to the role of nodes in non-fully-connected subgraphs. We then compare our approach to state-of-the-art models that can, in principle, capture the system's expected behaviour. Where fair comparisons are possible we show that our model displays excellent agreement with existing models, otherwise we show our model to either outperform existing models or to produce accurate results where other models fail. Finally, we use the generalised model to investigate the effect of loops/cycles as well as the impact of higher-order stucture, where global clustering is kept constant, on epidemic dynamics.

\section{Materials and Methods}\label{sec:mandm}

In this section we consolidate and generalise existing work centred around deriving low dimensional, deterministic and approximate ODEs that capture the time evolution of 
epidemic dynamics on configuration model networks. First, we re-introduce the basic susceptible-infected-recovered (SIR) epidemic model on random graphs following Volz's 
original PGF-based derivation \citep{volz2008sir}. This is followed by a rigorous formalisation of the hyperstub configuration model that was first presented by Karrer \& Newman \citep{karrerclust2010}. We then demonstrate how this model may be used to generate networks of differing subgraph compositions whilst keeping traditional network metrics such 
as first and second moments of the degree distribution, clustering and where possible the entire degree distribution, equal. Sec.~\ref{sec:oded} provides a derivation of the PGF-based approximate ODE model that accurately  captures SIR dynamics on hyperstub configuration networks. This derivation is similar to Volz et al.'s PGF-based extension from configuration and unclustered to clustered networks \citep{vmclust}, but generalised to incorporate arbitrary subgraphs. Finally, sec.~\ref{sec:TMA} provides an algorithm that automatically generates and solves ODEs presented in sec.~\ref{sec:oded} for SIR epidemics on networks constructed using a user-specified set of subgraphs.

\subsection{SIR epidemics on random graphs}\label{sec:SIRRND}

The SIR compartmental model involves a population with three types of individuals -- susceptible, infected or recovered -- whose interactions are modelled by a network. Infection travels across edges at a per-edge rate of $\tau$ and individuals recover, independently, at rate $\gamma$. To account for the heterogeneous contact patterns, the model is centred around the PGF induced by the network's degree distribution,
\begin{eqnarray*}
\psi(x) = \sum_{k = 0}^\infty p(k)x^k,
\end{eqnarray*}
where $p(k)$ is the probability that a randomly chosen node has $k$ links. Before we can demonstrate the usefulness of storing the network in this compact way, we need to define the \emph{survivor function}, $\theta(t)$. First, we define \emph{infectious contact} to be the event whereby an infected node $v$ transmits to its neighbour $u$, regardless of its state, i.e., irrespective of whether or not it is susceptible \citep{miller2011note}. Next, we select an edge uniformly at random, with nodes $u$ and $v$ at its ends, and define a direction from node $v$ to node $u$.
Let $\theta(t)$ be the probability that there has never been infectious contact from node $v$ to node $u$ by time $t$. Since an infectious contact does not depend on the state of the receiving node, Volz proposed the following simplifying assumption, ``we disallow infectious contact from node $u$ to node $v$''. Otherwise, $u$ may be infected by some other source, and in turn, infect $v$, thus increasing the probability of infectious contact from $v$ to $u$. This definition effectively implies that $\theta(t)$ is independent across all edges. For example the probability that a degree two node is susceptible at time $t$ is given by $\theta(t)^2$, or more generally
\begin{eqnarray*}
\psi(\theta(t)) =  \sum_{k = 0}^\infty p(k)\theta(t)^k =: S(t),
\end{eqnarray*}
where $S(t)$ is the fraction of susceptibles at time $t$. To analytically describe $\theta(t)$, we need to consider the rate at which a node with degree one becomes infected. This yields
\begin{eqnarray*}
\frac{d}{dt}\left(1 - \theta(t)\right) = \tau \theta(t) \frac{M_{SI}(t)}{M_{S}(t)} \Rightarrow \frac{d\theta(t)}{dt} =  -\tau \theta(t) \frac{M_{SI}(t)}{M_{S}(t)},
\end{eqnarray*}
where $M_S(t)$ and $M_{SI}(t)$ denote the expected degree of a susceptible node and the expected number of $SI$ edges per node at time $t$. Hence, $M_{SI}(t)/M_{S}(t)$ denotes the probability that a susceptible and infected node are connected at time $t$. In other words, a node which up to time $t$ is susceptible will, on average, become infected at rate $\tau M_{SI}(t)/M_{S}(t)$. It turns out that $M_S(t)$ can be computed using the PGF and is given by
\begin{eqnarray*}
\theta(t) \left.\frac{d\psi(x)}{dx}\right|_{\theta(t)} =  \sum_{k = 0}^\infty kp(k)\theta(t)^{k},
\end{eqnarray*}
which can be interpreted as the expected degree conditional on nodes being susceptible. To compute $M_{SI}(t)$ additional information from the PGF must be extracted, namely the \emph{excess degree}. This involves selecting an edge at random and following it to its originating node. The observed degree of this node, excluding the edge by which it was selected, is known as the \emph{excess degree} and has a distribution that is generated by
\begin{eqnarray*}
g(z) = \frac{1}{\langle k \rangle} \sum_{k = 0}^\infty (k+1)p(k+1)z^{k}.
\end{eqnarray*}
As before it is possible to condition this on susceptible nodes and thus to compute the expected excess degree of susceptible nodes
\begin{eqnarray*}
\theta(t) \left.\frac{dg(z)}{dz}\right|_{\theta(t)} = \frac{1}{\langle k \rangle} \sum_{k = 0}^\infty k(k+1)p(k+1)\theta(t)^{k}=: \delta_S(t).
\end{eqnarray*}
By assuming that the expected degree of a newly infected node is equal to the expected degree of a susceptible node, Volz uses the above, multiplied by $\tau$, to model the expected number of edges the disease can spread across upon infection of a susceptible node. This can be used to derive the equations that describe the flux between edges in different states. Namely, these are given by
\begin{eqnarray*}
\frac{d M_{SS}(t)}{dt} &=& - 2\delta_S M_{SS}(t), \nonumber \\ 
\frac{d M_{SI}(t)}{dt} &=& -M_{SI}(t)(\tau + \gamma) + 2\delta_S(t)M_{SS}(t) - \delta_S(t) M_{SI}(t), \nonumber 
\end{eqnarray*} 
where $M_{SI}(t)(\tau + \gamma)$, $2\delta_S(t)M_{SS}(t)$ and $\delta_S(t)M_{SI}(t)$ denote the $I$ infecting the $S$ or the $I$ recovering, $M_{SI}$ being created by a node in a $SS$ edge being infected by an external source to that $SS$ edge and, finally, the susceptible in a $SI$ edge being infected by an external source, respectively.
Summarising all the above yields the complete system of equations,
\begin{eqnarray*}
\frac{dS(t)}{dt}&=& \frac{d \theta(t) }{dt} \psi(\theta(t)), \\
\frac{dI(t)}{dt} &=& - \frac{d \theta(t) }{dt} \psi(\theta(t))- \gamma I(t), \nonumber\\
\frac{d M_{SS}(t)}{dt} &=& - 2\delta_S(t) M_{SS}(t), \nonumber \\ 
\frac{d M_{SI}(t)}{dt} &=& -M_{SI}(t)(\tau + \gamma) + 2\delta_SM_{SS}(t) - \delta_S(t) M_{SI}(t), \nonumber \\
 \frac{d\theta(t)}{dt} &=&  -\tau \theta(t) \frac{M_{SI}(t)}{M_{S}(t)}, \nonumber \\
 R(t) &=& 1 - S(t) - I(t). 
\end{eqnarray*} 
This concludes the derivation for PGF-based epidemic dynamics on random networks.  Volz et al. extended this methodology to clustered networks by defining a joint probability distribution which describes the typical number of lines and triangles allocated to nodes \citep{vmclust}. This particular derivation has been omitted from this paper. However, in the following section, we will outline a further generalisation of this whereby the joint probability specifies the distribution of subgraphs of various types around nodes. This then leads to more complex PGFs.

In App. \ref{sec:equiv}, we show how the PGF used in the main result of this paper can be made equivalent to the PGF resulting from Volz et al.'s original edge-triangle model.

\subsection{Hyperstub configuration model}\label{sec:hscm}

In this paper we generalise the configuration model~\citep{bollobas1980probabilistic} to the \emph{hyperstub configuration model}. Before we specify the model we need to establish how to classify hyperstubs, the set of stubs that connect a node to a subgraph, depending on their parent subgraph and their role within that subgraph.

To generate a hyperstub configuration network model one needs to first decide on a set of subgraphs or building blocks that will form the network. This is then followed by the identification of the number of different hyperstubs induced by the subgraphs: hyperstubs must be uniquely associated with both their parent subgraph and the \emph{orbit} of their incident nodes ~\citep{karrerclust2010} where the \emph{orbit} of a node is the set nodes with which it may be permuted such that no edges are created or destroyed. For example, in Fig.~\ref{fig:enu_hyperstub}, subgraph $G_\boxslash$ contains two distinct orbits $\{x_{14},x_{17} \}$ and $\{x_{15},x_{16} \}$.

Once all hyperstubs have been identified it is possible to define a joint  probability distribution that specifies the probability of a node having a certain combination of these. For example $f(x,y) = p_{x,y}$ may denote the probability of a node having $x \times G_0$ and $y \times G_\triangle$. Using this distribution it is possible to generate hyperstub degree sequences. For network generation these sequences will be subject to cardinality constraints. For example, the sum of the degree sequence of $G_\triangle$ must be divisible by three. Otherwise, the sequence needs to be re-generated. For asymmetric subgraphs, e.g.,  $G_\boxslash$, the sum of the degree sequences of both types of hyperedge must also be equal. In practice, this can be achieved by generating a suitable degree sequence for one type of hyperedge and then randomly permute it to obtain a second sequence for the second hyperstub. $G_\boxslash$ has two degree sequences, one for each hyperstub, and both must be even since we select pairs of nodes from each to form the subgraph.

The network generating algorithm will then form a dynamic list for each hyperstub, where a node with hyperstub degree $k_i$ appears $k_i$ times. This is followed by 
selecting nodes from the lists, at random and without replacement, and by following the subgraphs' hyperstub composition in order to construct subgraphs and the network. It is possible that self or multi-edges form in which case the selection is discarded and new samples chosen until a valid selection is obtained. This is repeated until all lists are empty.

In this paper we wish to both computationally generate networks and theoretically describe dynamics on such networks. The PGF of the hyperstub degree distribution provides the link between theory and simulation. The construction of the PGF induced by the hyperstub distribution can be achieved by encoding different levels of detail. At the simplest level nodes may belong to a number of subgraphs without further specifying their orbit or position within the subgraph \citep{vmclust}. The PGF could be constructed at the level of hyperstubs but would not differentiate between topologically equivalent positions in the subgraph, and this is what we use in our network generating algorithm (nodes may now be allocated asymmetric subgraphs) \citep{karrerclust2010}. Finally, the PGF can be specified by accounting for all details described above with the addition of the precise position of nodes within the subgraph (used in the ODE derivation, sec.~\ref{sec:oded}). For network generation the PGF takes the general form,
\begin{eqnarray}
\psi(\hat{z}) &=& \sum_{\hat{h}=0}^{\infty} p_{\hat{h}} \prod_{i=1}^m z_i^{h_i} \nonumber,
\end{eqnarray}
where $\hat{z}=(z_1,z_2,\dots,z_m)$ is a placeholder and $\hat{h} = (h_1,h_2,\dots,h_m)$ denotes the number of $h_i$ hyperstubs assigned to a node. The symbolic form of the PGF provides more flexibility for computation. Let us consider subgraphs distributed as follows: $G_0 \sim Pois(\lambda_1)$, $G_\triangle \sim Pois(\lambda_2)$ and $G_\boxslash \sim Pois(\lambda_3)$ (both hyperstubs of $G_\boxslash$ are Poisson distributed with parameter $\lambda_3$). The PGF of such a network is
\begin{eqnarray}
\psi(z_1,z_2,z_3) = exp\left(\lambda_1(z_1-1) + \lambda_2(z_2 - 1) +\lambda_3(z_3-1) \right) \nonumber.
\end{eqnarray}
From this PGF, the average number of subgraphs a node belongs to may be computed
\begin{eqnarray}
\left.\frac{\partial \psi(\hat{z})}{\partial z_1}\right|_{\hat{z}=1} = \lambda_1 =: \langle G_0 \rangle \nonumber.
\end{eqnarray}
By replacing $z_i$ with $z^a$, where $a$ is the number of stubs contained within the hyperstub $h_i$, the PGF of the classical degree distribution can be recovered
\begin{eqnarray}
\psi(z) := exp\left((\lambda_1(z-1) + \lambda_2(z^2-1) +\lambda_3(z^{5/2}-1) \right) \nonumber.
\end{eqnarray}
The $z^{5/2}$ term accounts for the fact that $G_\boxslash$ is counted twice, once for each of its hyperstubs. The first and second moments of the degree distribution are directly computed using the linearity of expectation and the fact that $Var(aX) = a^2X$. As well as recovering the degree distribution, it is possible to determine the expected number of triangles per node: $\langle \triangle \rangle  = \lambda_2 + 3/2 \lambda_3$, since on average each node in $G_\boxslash$ is incident to $3/2$ triangles. To summarise, we have
\begin{eqnarray}\label{eq:sys1}
\langle k \rangle &=&  \lambda_1 + 2\lambda_2 + \frac{5}{2}\lambda_3, \nonumber\\
Var(k) &=& \lambda_1 + 4 \lambda_2 + 25/4 \lambda_3,  \\
\langle \triangle \rangle  &=& \lambda_2 + 3/2 \lambda_3. \nonumber
\end{eqnarray}
By including a fourth subgraph in the above example,  the equivalent of system Eq.~\eqref{eq:sys1} will be underdetermined with 3 equations and 4 unknowns. 
This allows the first and second moments and the expected number of triangles (and therefore clustering) to be fixed whilst varying the subgraph composition. 
For example, fixing $\langle k \rangle =4$, $Var(k) = 8$ and  $\langle \triangle \rangle = 2$, we can form the underdetermined system,
\begin{eqnarray}
\begin{pmatrix}
 1 & 2 &  2 & 5\\
 1 & 4 & 4  & 25\\
 0 & 1 &  0 & 10
 \end{pmatrix}
\begin{pmatrix}
 G_0 \\
 G_\triangle \\
 G_\boxempty  \\
 G_{6c}
 \end{pmatrix}  
 =
 \begin{pmatrix}
 4 \\
 8 \\
 2  
 \end{pmatrix} , \nonumber
\end{eqnarray}
where the columns of the LHS matrix correspond to contributions to  $\langle k \rangle $, $Var(k)$ and  $\langle \triangle \rangle$ respectively and $G_{ic}$ denotes a complete subgraph of $i$ nodes. From this system it is possible to obtain two valid solutions: (1) $G_\triangle \sim Pois(2)$ and (2) $G_0 \sim Pois(9/2)$, $G_{6c} \sim Pois(3/10)$. Moreover, by replacing $G_{6c}$ with other types of subgraph and updating the L.H.S matrix, several differing network models with the same first and second moments and clustering may be obtained. A selection of such networks used in the results section is listed below:
\begin{eqnarray}\label{eq:model1}
\mbox{Model 1}:& G_\triangle \sim Pois(2), \nonumber\\
\mbox{Model 2}:& G_0 \sim Pois(2),~ G_\boxtimes \sim Pois(2/3), \nonumber \\
\mbox{Model 3}:& G_0 \sim Pois(8/3),~ G_{5c} \sim Pois(1/3), \nonumber \\
\mbox{Model 4}:& G_0 \sim Pois(3),~ G_{6c} \sim Pois(1/5). \nonumber
\end{eqnarray}

While the three most basic network metrics for the networks above are identical, their degree distributions are not. However, it is also possible to generate classes of networks where the degree distribution is equal between networks but the subgraph composition is not. Let us consider networks constructed purely out of cycles, where, regardless of the length of the cycle, cycle hyperstubs are composed of only pairs of stubs. It is then possible to increase the size of cycles whilst maintaining identical classical degree distributions  between different networks. This is implemented in the following way: first, allocate to each node, on average, a pair of cycle hyperstubs, then for each type of network allow the hyperstubs to form increasingly large cycles, starting with $G_\triangle$ then $G_\boxempty$ and so on. If the hyperstubs are distributed such that $h_i \sim Pois(2)$ then the classical degree distribution for each network will be such that only even degrees are possible, i.e., $P(degree=2k)=P(degree=k|Pois(2))$ denoted  $G_0 \sim 2 Pois(2)$ for convenience. It is also possible to include a null, random, model for comparison, i.e., a network with degree distribution given by $G_0 \sim 2 Pois(2)$ but connected at random. In our investigation we shall be using the following cycle-based networks:
\begin{eqnarray}\label{model2}
\mbox{Null Model}:& G_0 \sim 2Pois(2), \nonumber \\
\mbox{Model C1}:& G_\triangle \sim Pois(2), \nonumber \\
\mbox{Model C2}:& G_\boxempty \sim Pois(2), \nonumber \\
\mbox{Model C3}:& G_{\pentagon} \sim Pois(2), \nonumber \\
\mbox{Model C4}:& G_{\hexagon} \sim Pois(2), \nonumber
\end{eqnarray}
where $G_{\pentagon}$ and $G_{\hexagon}$ denote cycles of 5 and 6 nodes (pentagons and hexagons), respectively. 
Having thus created two classes of networks, the former will be used to show how conventional network metrics may not entirely capture the structure of the network as far as dynamics are concerned; the latter to investigate the effect of cycles of increasing length on dynamics.

\subsection{SIR epidemics on hyperstub configuration model networks}\label{sec:oded}

This section presents the derivation of a general $SIR$ epidemic model for a network built from an arbitrary number of subgraph types. Conceptually, this model uses the node labelling approach of \citep{karrerclust2010} and generalises the PGF-type framework of Volz et al. \citep{vmclust,volz2008sir}. By taking this approach it is possible to derive ODEs that accurately predict the epidemic prevalence on networks that exhibit a variety of exotic subgraphs, both fully- and non-fully connected. 

\begin{figure}[h]
\centering
 \includegraphics[scale=0.95]{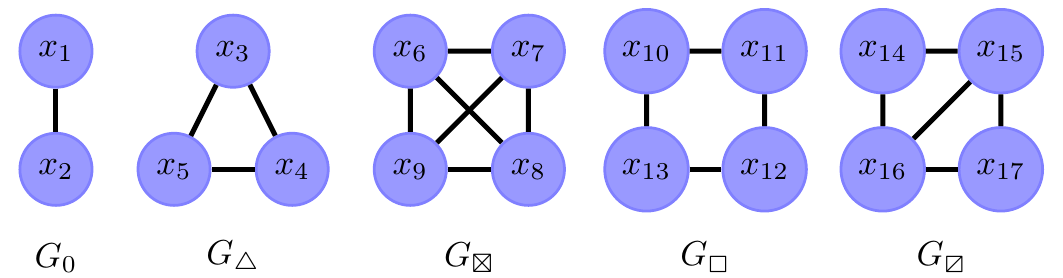} 
\caption{Subgraph notation and position labelling. Subgraphs are denoted by $G$ followed by a symbolic subscript for ease of reference. }
 \label{fig:enu_hyperstub}
\end{figure}

The first step is to choose the set of subgraphs to be included in the network. Let an arbitrary set of subgraphs be labelled by $\{G_1, G_2, \dots, G_{M}\}$. For example, Fig.~\ref{fig:enu_hyperstub} shows $M=5$ different subgraphs, which result in $m=17$ distinct node positions, where $m$ stands for the total number of nodes over all subgraphs. For clarity, we recall that a hyperstub is the set of half-links connecting a node to a subgraph. This example highlights the key component of the model, namely to distinguish between all nodes of a subgraph even those that are topologically equivalent. This distinction makes it possible to deal with the added complexity of having to account for labelled subgraphs. Each node/position of a subgraph is labelled. This is reflected in a PGF that accounts for each and every node in each and every subgraph. This gives rise to a PGF of the following form:
\begin{eqnarray}\label{eq:hpgf}
\psi(\hat{\alpha}) &=& \sum_{\hat{y}=0}^{\infty} p_{\hat{y}} \prod_{i=1}^m \alpha_i^{y_i}, \nonumber
\end{eqnarray}
where $\hat{\alpha}=(\alpha_1,\alpha_2,\dots,\alpha_m)$ is a placeholder and $\hat{y} = (y_1,y_2,\dots,y_m)$ is such that $y_i$ is the number of times a node appears in position $x_i$, $i=1, \dots, m$. 

For each subgraph its state at time $t$ is denoted by $G_{x}(S,I,\dots,R)$. This not only describes a subgraph and its state but also the \emph{expected} number of the given subgraph in the given state at time $t$, i.e., when appended with a state this notation has numerical meaning. Since $G_{x}(S,I,\dots,R)$ accounts for the state of node, it will always explicitly depends on $t$. To describe the flux between different subgraph states, infectious events within and \emph{between} subgraphs need to be considered. This requires a generalisation of $\theta(t)$ which was first given in sec.~\ref{sec:SIRRND}. Accordingly, we now first select a hyperstub at random and then define a direction, from its parent subgraph to its incident node. An \emph{infectious contact} is now the event that $u$, regardless of its state, becomes infected by one of its adjacent nodes within that subgraph. $\theta(t)$ now needs to reflect a node's position in the subgraph. Hence, we define  define $\theta_i(t)$ to be the probability that the group of edges connecting a node $u$ in position $x_i$ to the parent subgraph have not allowed for \emph{infectious contact} from any infectious node in the subgraph to $u$ by time $t$. Again, we impose that $u$ cannot transmit infection to the subgraph in question. Under these assumptions, the infectious contact through hyperstubs to position $x_i$ is now independent. A node that appears only $k$ times in position $x_i$ remains susceptible with probability $\theta_i^{k}(t)$. By geometrically compounding all $\theta_i(t)$ into a PGF, it is possible compute the fraction of the susceptible population. This is given by
\begin{eqnarray}\label{eq:susceptible1}
S(t) = \psi(\hat{\theta}) = \sum_{\hat{y}=0}^{\infty} p_{\hat{y}} \prod_{i=1}^m\theta_i^{y_i}.
\end{eqnarray}
This probability is equal to the fraction of susceptible nodes in the population at time $t$ \citep{volz2008sir}. $\theta(t)$ is referred to as a survivor function. It is dependant on time and may by computed from first principles using the definition of the Poisson process. However, in our formulation, it is computed from variables that denote the expected rate, $T_i$, at which infection is transmitted to a node in position $x_i$ through the corresponding subgraph. We note that while $T$ is commonly used to denote the cumulative probability that infection may occur, we keep it as defined above to be consistent with the current literature on such models \citep{vmclust}. Each position label $x_i$ has a $T_i$ variable associated with it. The following examples show these rates for positions $x_1$, $x_2$ and $x_3$, see Fig.~\ref{fig:enu_hyperstub}:
\begin{eqnarray}\label{eq:rofg}
T_1 &=& \tau[G_0(SI)],\\
T_2 &=& \tau[G_0(IS)],\\
T_3 &=& \tau[G_\triangle(SSI) + G_\triangle(SIS) + 2G_\triangle(SII) \nonumber \\
    & & +G_\triangle(SRI) + G_\triangle(SIR)].
\end{eqnarray} 
To generate the above identities, we consider a susceptible node in position $x_i$ and list all possible corresponding subgraph states that allow this node to be exposed to infection. 
$T=(T_1, T_2, \dots, T_{m})$ can now be used to determine the probability that a susceptible node has an infectious neighbour within a certain subgraph type. 
This is done by dividing $T_i\tau^{-1}$ by the number of states that involve a susceptible at position $x_i$:
\begin{eqnarray}\label{eq:adegree}
\frac{T_i}{\tau \sum\limits_{A,B,C,D} G_{(\cdot)}(x_i=S,\dots,A,B,C,D)}. \nonumber
\end{eqnarray}
The expected degree of a susceptible node at position $x_i$ is given by
\begin{eqnarray}
\langle k_i \rangle = \sum_{\hat{y}=0}^{\infty} y_i p_{\hat{y}} \prod_{i=1}^m \theta_i^{y_i} = \theta_i \left.\frac{\partial \psi}{\partial \alpha_i}\right|_{\alpha = \hat{\theta}}, \nonumber
\end{eqnarray}
where $\hat{\theta} = (\theta_1, \theta_2, ... , \theta_m)$. To compute the expected degree for every position of every subgraph, one can take the Jacobian of $\psi$ evaluated at $x = \hat{\theta}$,
\begin{eqnarray}
J(\psi)|_{\alpha=\hat{\theta}}. \nonumber
\end{eqnarray}
The $i^{th}$ entry of this vector evaluated at $\alpha=\hat{\theta}$ shall be denoted $J_i$. A susceptible node in position $x_i$ will have remained susceptible up to time $t$, with probability $\theta_i$ after which infection may be transmitted at rate $T_i/J_i$. This information may be used to form the following equation: 
\begin{equation}\label{eq:theta}
\frac{d}{dt}(1-\theta_i(t)) = \theta_i(t) \frac{T_i}{J_i} \Rightarrow \frac{d \theta_i(t)}{dt} = -\theta_i(t) \frac{T_i}{J_i}.
\end{equation}
$\dot{\theta}_i(t)$ decays at the rate at which a subgraph transmits infection to its node in position $x_i$, conditional on that node being susceptible. 

\begin{figure}[h]
\centering
 \includegraphics[scale=0.95]{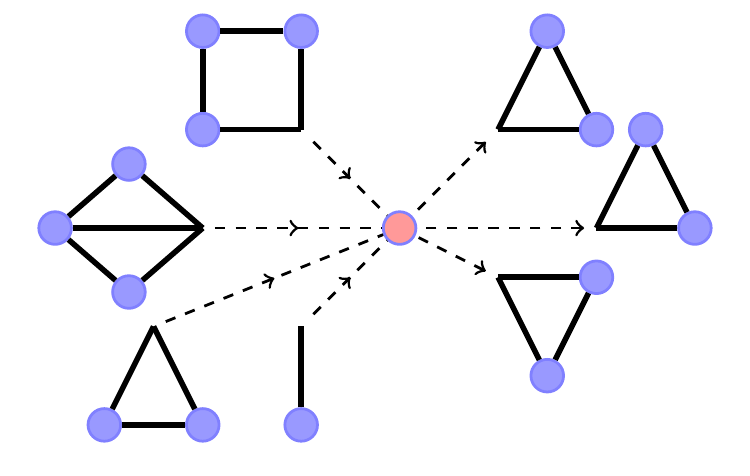}
\caption{Graphical representation of $(T\Delta)_i$. $\Delta$ and $T$ denote the excess degree of a susceptible node and rate of infection, respectively. We note that newly infected nodes are modelled as previously susceptible nodes so the product $(T\Delta)_i$ is being used to model the expected number of $x_i$ edges infection will be able to spread along upon infecting a susceptible node. This product implicitly considers all possible routes of infection into the node. The left hand side of the figure shows example subgraphs that are the source of infection for the red node. The right hand side of the figure graphically represents the expected excess degree of $G_\triangle$ subgraphs for the red node.}
\label{fig:flux}
\end{figure}

Once a node is newly infected it is important to determine what, if any, subgraph states are created or destroyed. To do this, we use the susceptible nodes's excess degree prior to the infection. For the full derivation of susceptibles' excess degree refer to App.~\ref{sec:excess}. In this derivation, the excess degree must be generalised to account for the degree of the different positions a node may be in, i.e., $\langle k_i \rangle,~ i=1,2,\dots,m$. The expected excess degree for susceptible nodes is given by
\begin{eqnarray}
\Delta_{i,j} = \theta_j \left.\frac{H_{i,j}(\psi)}{J_i(\psi)} \right|_{\alpha = \hat{\theta}}, \nonumber
\end{eqnarray}
where $H(\psi)$ is the Hessian of the PGF. $\Delta_{i,j}$ denotes the expected number of $x_j$ positions associated with a node that has been selected at random, but proportionally to the number of $x_i$ positions associated with that node. It is now possible to formulate ODEs describing the evolution subgraph states. We denote the time derivative of a subgraph's state by $\dot{G}_{(\cdot)}$. This quantity is dimensionless but not normalised. For example, the number of unique $(SI)$ links in a network of size $N$ is given by $[SI] = NG_0(SI)$. To form the ODE for the subgraph state $G_{0}(SI)$, we consider all possible ways in which this state may be created or destroyed, namely 
\begin{eqnarray}
\dot{G}_0(SI) &=& -(\tau + \gamma)G_0(SI) \label{eq:example1} \nonumber\\ 
			&& - (T \Delta)_1 G_0(SI) + (T \Delta)_2 G_0(SS), \label{eq:example2}
\end{eqnarray}
where $(T \Delta)_1$ denotes the first entry of the vector that is the product of the matrix $\Delta$ multiplied from the left by vector $T$. Conceptually $(T \Delta)_i$ denotes the expected number of nodes in position $x_i$ an infection will encounter upon infecting a susceptible node through any possible route, see Fig.~\ref{fig:flux}. The first term on the RHS of Eq. \eqref{eq:example1} describes this state being destroyed by the $I$ infecting the $S$ or the $I$ recovering. The second term stands for this state being destroyed by the $S$ being infected by an outside source. Finally, the last term corresponds to this state being created by the second node of $G_0(SS)$ being infected by a source external to the subgraph. To further illustrate this, the equations for $G_0(SS)$ and $G_0(IS)$ are given,
\begin{eqnarray}
\dot{G}_0(SS) &=& -[(T \Delta)_2+ (T \Delta)_1] G_0(SS), \nonumber\\
\dot{G}_0(IS) &=& -(\tau + \gamma)G_0(II) \nonumber \\
  			  & & -(T \Delta)_2 G_0(IS) +(T \Delta)_1 G_0(SS). \nonumber
\end{eqnarray}
Equations for every state of every subgraph must be derived. In general, we first describe any infection and recovery events of nodes within a subgraph. 
Next we list all possibilities for susceptible nodes to be infected from sources external to that subgraph using the appropriate $(T \Delta)$ terms. 

To compute network-level prevalences, we recall that $S(t)$ can be computed at any time by Eq.~\eqref{eq:susceptible1}. $\dot{I}(t)$ is computed directly by differentiating $S(t)$. Namely, since susceptibles become infected and and infected nodes recover at rate $\gamma$, we have
\begin{eqnarray}\label{eq:ep}
\dot{I}(t) &=& \sum_{i=1}^m \dot{\theta}_i(t) \frac{\partial \psi (t)}{\partial \theta_i} - \gamma I(t), \\\dot{R}(t) &=& \gamma I(t).
\end{eqnarray} 
The total number of equations is given by $2 + m + \sum_{i=1}^M 3^{|G_i|}$, where $|\cdot|$ denotes the number of nodes in a subgraph.
In App.~\ref{sec:exampleodes} we give or more example ODEs and in App.~\ref{sec:equiv} we show how our model is equivalent to 
previous systems developed for complete subgraphs \citep{vmclust}.
\subsection{Initial conditions} \label{sec:intial}

Let $\epsilon $ be the fraction of initially infected nodes. Hence, $\epsilon = I_0/N$, where $I_0$ is the number of initially infected nodes and $N$ is the network size. Initial conditions for the $I$ and $R$ populations are given by
\begin{eqnarray}
I(0) = \epsilon,~ R(0) = 0. \nonumber
\end{eqnarray}
At time $t=0$  no hyperstub has transmitted infection, therefore, $\theta_i(t=0) = 1$.
For a subgraph that contains a single infected node, $G(t=0) = \epsilon \langle k \rangle_i $ where $\langle k \rangle_i$ is the expected hyperstub degree.
For the subgraph with every node susceptible we set $G(t=0) = (1-\epsilon) \langle k \rangle_i $. By assuming that only a small fraction of the population, i.e., a single node, is initially infected, we do not allow non-zero initial conditions for subgraphs with more than one infectious node.

\section{Automated code-generation of the mean-field model}\label{sec:TMA}

We now present our methodology for computationally generating a complete system of equations for a network constructed from subgraphs following a configuration model. This procedure requires the PGF of a hyperstub degree distribution (HDD), the adjacency matrices of corresponding subgraphs, and epidemiological parameters as inputs. The algorithm will output the system of ODEs that will predict the network-level prevalence. Table \ref{table:vars} gives a brief summary of the variables that need to be generated, listed in the order they are generated in this section.

Let $\vec{G}$ denote the vector of states of a subgraph $G$ with $\vec{G}_i$ denoting a specific state of $G$. For the $SIR$ model, $\vec{G}$ has $3^{|G|}$ elements. To generate $T_i$ from $\vec{G}$, the following steps are needed: (1) cycle through $\vec{G}$, (2) for each infectious contact to node $i$ in state $\vec{G}_j$, update $T_i$ to $T_i = T_i + \vec{G}$. Using $T$ the survivor functions can be computed, see Eq.~\eqref{eq:theta}, which are then used to compute the fraction of the population which is susceptible, infected or recovered, see Eq.~\eqref{eq:ep}.

The ODEs corresponding to subgraphs  need to be represented with a \emph{rate matrix}, $\mathbf{Z}$. This matrix encodes all information relating to the given subgraph, namely the excess degrees, rates of infection over subgraphs $T$, epidemiological parameters $\tau$ and $\gamma$, and implicitly encodes the subgraph's adjacency matrix $g$. To compute $\Delta$, we use Eq.~\eqref{eq:hpgf} and a symbolic software package to calculate the Jacobian and Hessian of the PGF.

For each subgraph, we initialise the matrix $\mathbf{Z}$ as a square matrix with all entries set to zero. The $i^{th}$ column and row of $\mathbf{Z}$ correspond to state $\vec{G}_i$. Once populated, the entry $\mathbf{Z}_{i,j}$ contains the rate at which state $i$ transitions to state $j$.

To illustrate how to generate $\mathbf{Z}$, we consider the $G_0$ subgraph, see Fig.~\ref{fig:enu_hyperstub}, with states $\vec{G} = (SS,$ $SI,$ $SR,$ $IS,$ $II,$ $IR,$ $RS,$ $RI,$ $RR)$. We associate the state $G_0(SS)$ with the first row and column of $\mathbf{Z}$. Moving along the top row, when a column index is reached that corresponds to a state that $G_0(SS)$ may transition to, we update the entry with the appropriate rate. The first row of $\mathbf{Z}$ is all zero except for $\mathbf{Z}_{1,2}=(T\Delta)_2$ and $\mathbf{Z}_{1,4}= (T\Delta)_1$. The second row, corresponding to state $G_0(SI)$, has entries $\mathbf{Z}_{2,3}=\gamma$ and $\mathbf{Z}_{2,5}=\tau + (T\Delta)_1$, see Eq.~\eqref{eq:example2}. Fill every row of the matrix $\mathbf{Z}$ in this way, refer to the SI text for the full matrix corresponding to $G_0$. The algorithm for this process is given for an arbitrary subgraph in App.~\ref{sec:alg2}, and the corresponding Matlab code is provided as supplemental material.

Using the rate matrix, the ODE for the subgraph state $\vec{G}_i$ yields
\begin{eqnarray}\label{eq:eq_gen}
\frac{d \vec{G}_i}{dt} = - \left(\sum_{j=1}^{3^{|G|}} \mathbf{Z_{i,j}} \right)\vec{G}_i + \left(\sum_{k=1}^{3^{|G|}} \mathbf{Z_{k,i}}\right)\vec{G}_k.
\end{eqnarray} 
The final step to generating the full system is to set the initial conditions. 
Only the initial conditions for subgraph states need computing as $I(0),~R(0)$ and $\theta_{i}(0)$ are fixed as per the previous section.
This can be done by cycling through each element of $\vec{G}$. If (a) $\vec{G}_i$ is a 
purely susceptible state then we set $\vec{G}_{i_{0}}=J_i(1-\epsilon)$, and if (b) $\vec{G}_i$ contains a single infectious individual and is otherwise susceptible, we set $\vec{G}_{i_{0}}=J_i \epsilon$. All other states are set to zero, as we assume that with a sufficiently small infectious seed, the probability of having two infectious individuals in a subgraph is zero.
\section{Results}

To validate the proposed mean-field model and to assess the goodness of the approximation, we compare results from the ODEs to output from stochastic simulations. 
Networks were generated following the configuration algorithm, please refer to App.~\ref{sec:alg1}. Typically we generated 500 networks of size $N=15000$ and computed a single realisation of the epidemic, according to the Gillespie algorithm with the per link rate of infection $\tau = 1$ and a recovery rate of $\gamma = 1$. 
Simulations which died out before an outbreak occurred were removed. The simulations were seeded with a single infectious individual and an outbreak was said to occur if 5\% infectious prevalence was achieved. In all plots simulation results and the solution of ODEs are plotted in solid lines and discrete points, respectively. 

\begin{figure}[ht!]
\begin{center}
 \includegraphics[scale=0.6]{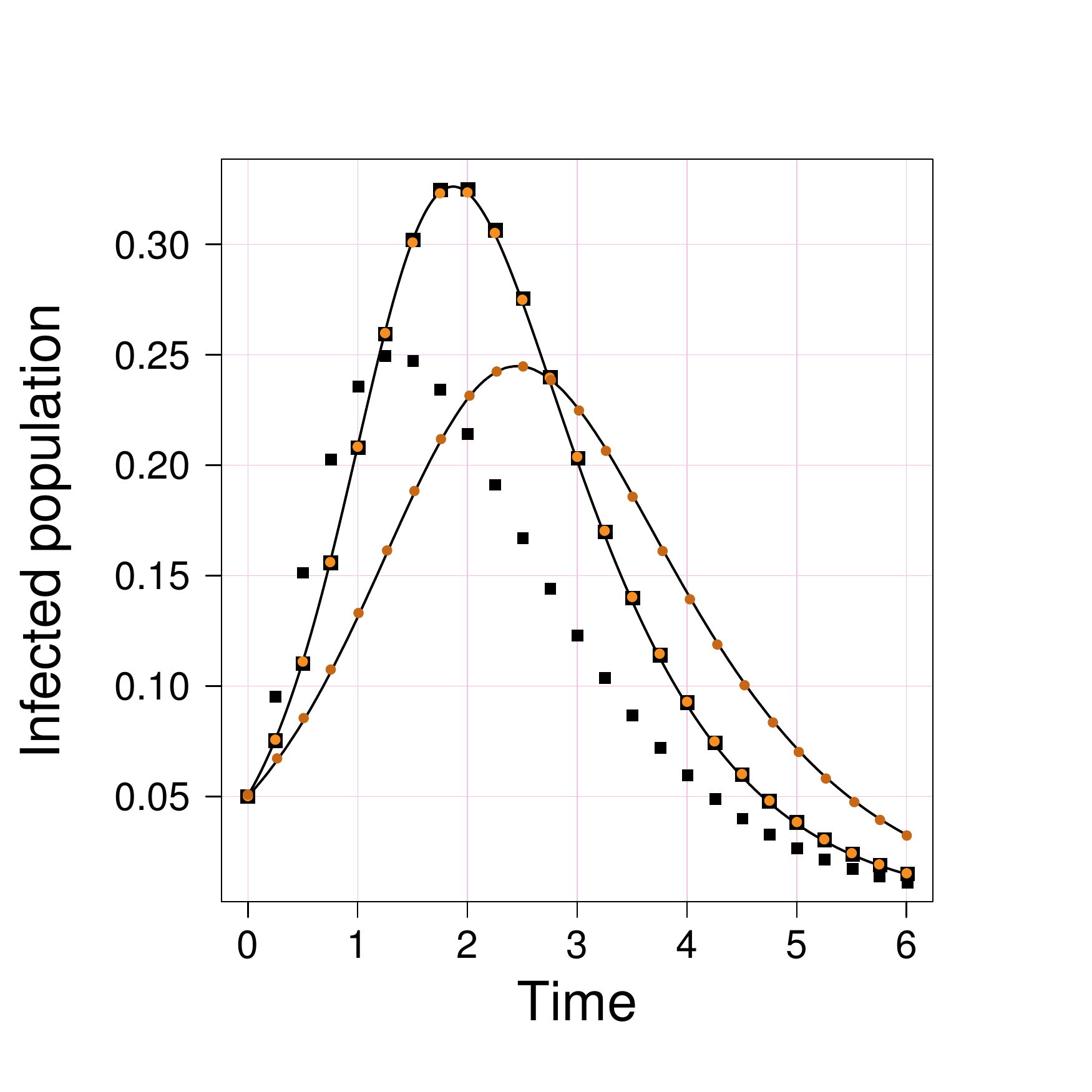}
 \end{center}
 \caption{Performance of other models. Lines, circles and squares correspond to simulation average, ODE solution and pairwise ODE solution, respectively. All networks are homogeneous with $k=5$. The lower peaks correspond to networks generated with each node allocated one of each corner type of a $G_\boxslash$ with clustering $\phi = 0.3$. Data with higher peak correspond to networks generated with a single $G_0$ and two $G_\boxempty$ subgraphs yielding $\phi \approx 0$.}
 \label{fig:allplotsa}
 \end{figure}

To start, we test the performance of our model against existing or state of the art models. To do this,
in Fig.~(\ref{fig:allplotsa}), we show results for two degree distributions that are homogeneous in the classical sense. Their PGFs are given by
\begin{eqnarray}
\psi_1(\hat{\alpha}) &=& \frac{1}{2}(\alpha_{14}+\alpha_{17}) \frac{1}{2}(\alpha_{15}+\alpha_{16}), \nonumber\\
\psi_2(\hat{\alpha}) &=& \frac{1}{2}(\alpha_1+\alpha_2) \frac{1}{4^2}(\alpha_{10}+\alpha_{11}+\alpha_{12}+\alpha_{13})^2,\nonumber
\end{eqnarray}
where the variables $\alpha_i$ correspond to subgraphs given in Fig.~\ref{fig:enu_hyperstub}. Figure~\ref{fig:allplotsa} shows results from a pairwise model with closures at the level of quadruples \citep{house2009motif,wrap5586}. 
While the classical clustering is easy to compute, the order-four clustering/transitivity ratios were measured following a recently developed subgraph counting algorithm \citep{ritchie2014higher}. These are defined as the ratio of a given subgraph count to all open and closed paths of length four, both counted uniquely. Currently, this model operates using an average or homogenous degree and stores no information about the degree distribution, but does assume random mixing of subgraphs. 

All models perform well in capturing the epidemic dynamics on networks generated using the PGF given by $\psi_1$, see Fig.~(\ref{fig:allplotsa}) with higher epidemic peak. However, when networks are created using the PGF given by $\psi_2$, see lower peak in Fig.~(\ref{fig:allplotsa}), the pairwise model struggles to accurately capture the dynamics, both anticipating and compressing the epidemic's time scale or duration, but predicting correctly the peak prevalence. The pairwise model does not encode any information relating to degree or subgraph distribution and hence a homogeneous random set-up, as used here, would be an appropriate choice.
  
\begin{figure}[ht!]
\begin{center}
 \includegraphics[scale=0.6]{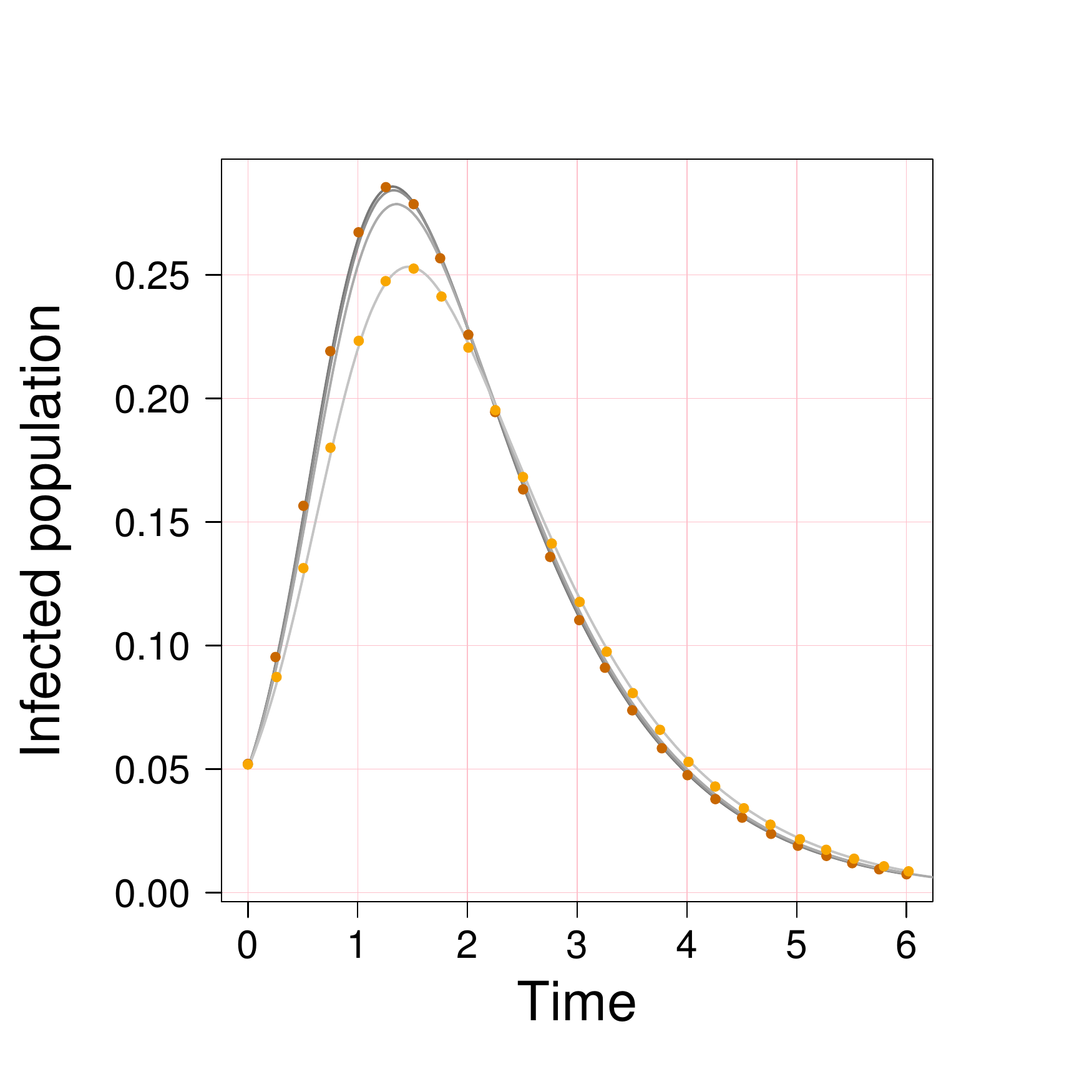}
 \end{center}
 \caption{Clustering and cycles. Solid lines and markers correspond to simulation average and ODE solution, respectively. From darkest to lightest, the solid lines correspond to: $k \sim 2Pois(2)$, $G_{\pentagon} \sim Pois(2)$, $G_{\boxempty} \sim Pois(2)$ and  $G_{\triangle} \sim Pois(2)$, i.e., each network used has an identical degree distribution given by $P(degree=2k)=P(degree=k|Pois(2))$. Clustering is $\phi = 0.2$ and $\phi \approx 0$ for the $G_\triangle$ 
and other networks, respectively. For clarity, ODE solutions for only the two extreme cases, the null and triangle network, have been included.
Epidemics corresponding to cycles of length six have been computed but omitted due to their close similarity to the null case.
}
 \label{fig:allplotsb}
 \end{figure}

The key advantage of our algorithm over existing ones is that it can handle non-fully connected subgraphs. To test this, in Fig.~\ref{fig:allplotsb}, we utilise networks models C1-C4 as specified in sec.~\ref{sec:hscm}. Fig.~\ref{fig:allplotsb} shows plots of simulation average compared to the ODE's solution for the four network types. 
We observe that the epidemic behaviour of networks composed of increasingly large cycles quickly converge to that of the random null case.
It has previously been observed that for networks with the same degree distribution, an increasing level of clustering slows the epidemic transmission and requires a higher transmission rate in order to observe a successfully spreading epidemic \citep{keeling1999effects,green2010large}. This occurs for two reasons: (1) subgraphs that are densely connected share fewer connections to the rest of the network so an initial seed will be restricted to one part of the network and (2) this same effect leads to infectious nodes competing for susceptible nodes. While this may make transmission more efficient locally, it does limit further seeding in fully susceptible parts of the network. Fig.~\ref{fig:allplotsb} shows that the effect of $G_\boxempty$ is similar to that of the clustered network, but less pronounced; both the time and size of peak infectious prevalence is delayed and reduced when compared to the null case. For cycles larger than four nodes this behaviour is less pronounced and the epidemics for larger cycles converge to the null case, as observed with $G_{\pentagon}$.

To highlight the flexibility of our model and its wide-ranging applicability to systematically investigating the impact of higher-order network structure, in Fig.~(\ref{fig:allplotsc}), we consider four networks with the same first and second moments, and the same classical clustering but generated using different families of subgraphs, see models 1-4 sec.~\ref{sec:hscm}. Figure~\ref{fig:allplotsc} shows simulation average for all four networks and the solution of ODEs for the upper and lower cases, models 1 and 4 respectively. 

Figure~\ref{fig:allplotsc} shows a clear trend whereby larger subgraphs lead to epidemics with smaller peak prevalence. A second more subtle trend shows a delay in time until peak prevalence. Subgraphs of larger size lead to a significant difference in the behaviour of epidemics and echo what was observed for increasing levels in clustering. This could be explained by considering a subgraph with average degree $\langle k_s \rangle$. When $ \langle k \rangle < \langle k_s \rangle$ the network will exhibit extreme clustering, where isolated structures are increasingly densely connected at the cost of becoming isolated. This effect is more subtle than clustering but it can be significant. This suggests that the accuracy of future models would improve if they can correctly account for networks' subgraph composition, particularly subgraphs beyond that of triangles.

Finally, the data in Fig.~\ref{fig:allplotsc} has been produced using networks that do not have the same degree distribution but do have equal first and second moments, and clustering. To better understand how the non-equal higher moments may have affected the results, we have simulated epidemics on the corresponding random networks, Model $1':~G_0 \sim 2Pois(2)$ and Model $4':~G_0 \sim Pois(3) + 5Pois(1/5)$, see App.~\ref{sec:null}. This plot shows that the differences observed in Fig.~\ref{fig:allplotsc} cannot be explained by the difference in the degree distribution alone. Thus, generating identical clustering but using different subgraphs can lead to non-negligible differences in epidemic dynamics.

\section{Discussion}

Higher-order structures, captured for example as different subgraph compositions and arrangements in a network, have been identified as features of real networks. Examples include households, social interactions and biological networks. These building blocks of networks have been shown to play a key role in defining a network's topology and can have significant impact on the functions of the network or on the dynamical processes unfolding on the network. Despite this, the modelling toolset in this direction is underdeveloped. Here, we provided an approach that considerably extends the scope of the current modelling framework by enabling us to consider arbitrary sets of exotic subgraphs as building blocks for the network. Our approach also offers control over the arrangements of subgraphs and, more importantly and uniquely, an automated way of generating a system of ODEs that accurately capture the prevalence profile for a wide range of subgraph sets, as shown in the results section.

The previous section has shown how higher-order structures may be investigated using this model. Moreover, we provided the first example of generating classes of networks constructed using different subgraph sets while keeping degree, variance and clustering, all in the classic sense, fixed. For example, we showed that epidemics on networks with no clustering, but exhibiting open loops, display features which are significantly different to those observed in classical random networks with effectively no clustering. Equally, we have shown that different subgraph combinations or arrangements can create higher-order structure that may significantly affect the epidemic dynamics. 
Our work opens the possibility to carry out a wide-ranging and systematic investigation of the impact of subgraphs and higher-order structure on dynamics on networks. When presented with real world network data whose structure can be explained by a set of subgraphs, all that will be needed in order to apply our framework is to extract the subgraphs and their distribution around nodes. A possible limitation to the widest applicability is the number of nodes in the largest subgraph. However, as shown by our results when going from squares to pentagons, it is likely that the effect of higher-order structures will decay, or be less marked, as their size increases.  

There are two key ways in which this work may be extended: (a) generalisation to $SIS$ dynamics. Due to the definition of $\theta(t)$ it is currently not possible to apply this model to $SIS$ dynamics. However, all the framework relating to network structure is independent from this variable and may therefore still be appropriate. (b) The subgraph approach is highly suitable for adaptation to household models. Household models typically specify a distribution of household sizes overlaid on a contact network to produce a well-connected network \citep{house2009household,ball2012sir}. A successful incorporation of such network in our framework could lead to a highly relevant set of household models.


\section{Figures and tables}

\newpage
\newpage

\newpage

\newpage 
 
 \begin{figure}[ht!]
\begin{center}
 \includegraphics[scale=0.6]{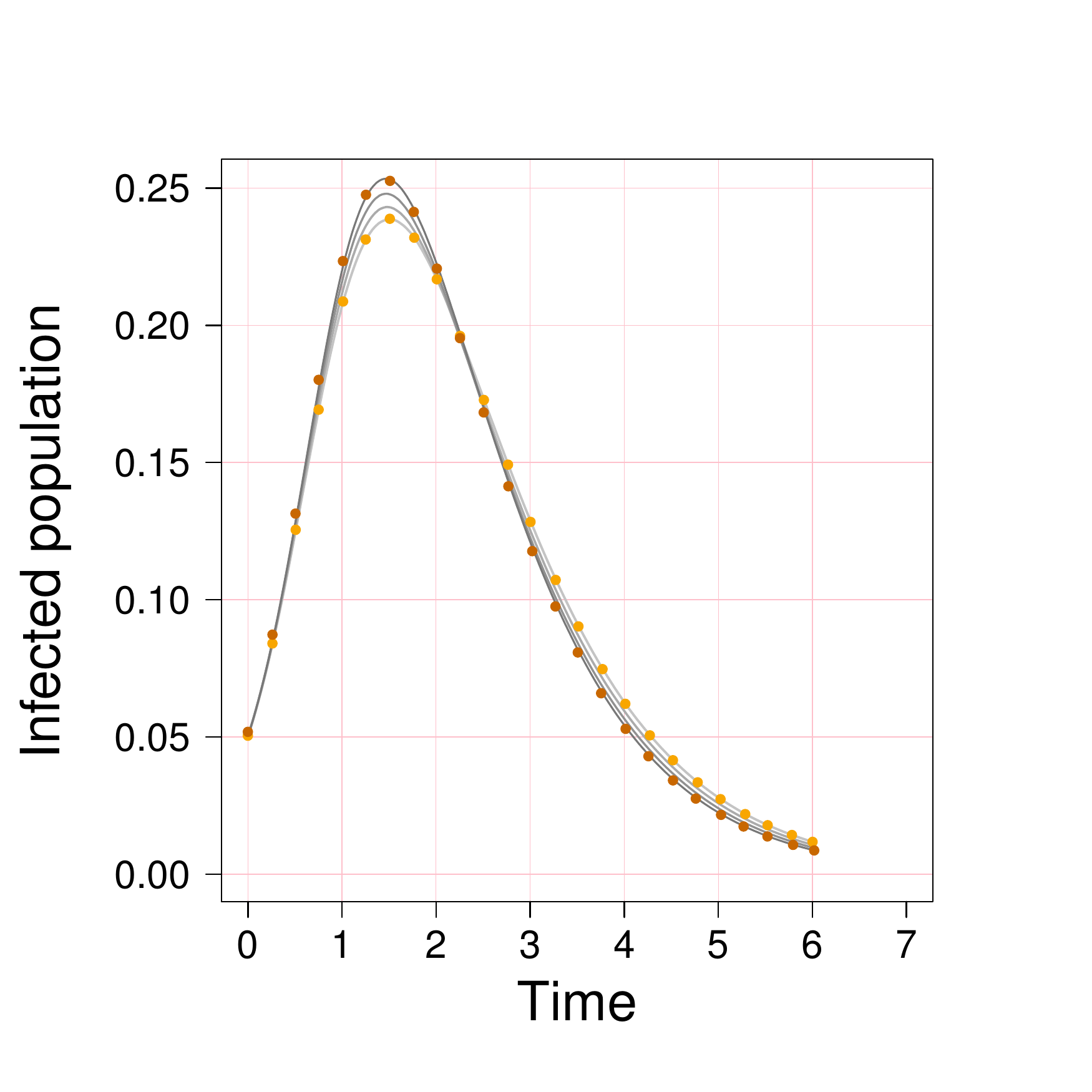}
 \end{center}
  \caption{Clustering via differing subgraphs. Solid lines and markers correspond to simulation average and ODE solution, respectively. From darkest to lightest the solid lines correspond to: $G_\triangle \sim Pois(2)$;  $G_0 \sim Pois(2)$,  $G_{\boxtimes} \sim Pois(2/3)$;  $G_0 \sim Pois(8/3)$, $G_{cp}, \sim Pois(1/3)$ and $G_0 \sim Pois(3)$, $G_{ch}, \sim Pois(1/5)$, where $cp$ and $ch$ denote complete pentagon and hexagon subgraphs, respectively. The networks were generated so that $\langle k \rangle = 4$, $var(k) = 8$ and $\phi = 0.2$. The downward trend of peak prevalence corresponds to networks composed of complete subgraphs of increasing size. 
 The larger subgraphs lead to more connections within the group rather than to the rest of the network.}
 \label{fig:allplotsc}
 \end{figure}

\newpage

\begin{table}[ht!]
\begin{center}
\begin{tabular}{p{2cm} p{3.8cm} p{3.8cm}}
Variable & Description & Generation \\  \hline 
\rule{0pt}{2.5ex}     $ \psi$ & PGF of the HDD given as a function, not as a series. & A symbolic software package can be used to compute the Jacobian and Hessian. \\ [0.5ex] 
 $\theta_{i}(t)$               & Survivor functions with their evolution equations given by ODEs.  & These ODEs can be defined within a single \texttt{for} loop, see Eq.\eqref{eq:theta}.  \\ [0.5ex] 
 $(S,I,R)$              & The prevalences of $S$, $I$ and $R$, with the latter two given by numerical solutions of ODEs & From Eq.~\eqref{eq:ep}, it follows that $S =\psi(\theta)$.    \\ [0.5ex] 
 $T_i$                  & Total rate of infection experienced by an $S$ in position $x_i$.  & For a subgraph with $m$ nodes, $T_i$ may be generated by $m$ nested \texttt{for} loops cycling through the possible states that a subgraph can be in, see Eq.~\eqref{eq:rofg}. \\ [0.5ex] 
 $G_x(S,I,\dots,R)$      & Expected prevalence of a subgraph in a given state.   & The equation for this is computed based on the rate matrix, $\mathbf{Z}$, see Eq.~\eqref{eq:eq_gen}.    \\[0.5ex] 
\end{tabular}
\label{table:vars} 
\end{center}
\caption{Summary of the key system variables and their generation.}
\end{table}

\newpage

\section{Appendix}

In this Appendix we (a) give a more detailed explanation of the excess degree, (b) provide ODEs for an example network, (c) show how our generalised model reduces to a previous model under certain conditions, (d) provide an example state transition matrix, (e) give pseudocode for both the subgraph-based configuration model and the algorithm used to generate the state transition matrix and, finally, (f) compare epidemic dynamics on two configuration model networks with their degree distributions being different but with the same mean and variance.

\subsection{Excess degree}\label{sec:excess}
Recall the probability generating function (PGF) of a network's hyperstub degree distribution with $m$ nodal positions:
\begin{eqnarray}\label{eq:Ahpgf}
\psi(\hat{\alpha}) &=& \sum_{\hat{y}=0}^{\infty} p_{\hat{y}} \prod_{i=1}^m \alpha_i^{y_i},
\end{eqnarray} 
where $\hat{\alpha}=(\alpha_1,\alpha_2,\dots,\alpha_m)$ is a placeholder, and $\hat{y} = (y_1,y_2,\dots,y_m)$ is such that $y_i$ denotes the number of times a node appears in position $x_i$, $i = 1,2, \dots,m$. The PGF of the excess degree distribution is a critical component in our derivation and it is illustrative to see how it is computed. To see this, we first compute the expected excess degree, select a node at random but proportional to its number of $x_i$ hyperstubs, $y_ip_{\hat{y}}$. Next, to obtain the expected $x_i$ degree, sum this product over all nodes 
\begin{eqnarray}
\langle x_i \rangle &=& \sum_{\hat{y}=0}^{\infty} y_i p_{\hat{y}} =: \left.\frac{\partial \psi(\hat{\alpha})}{\partial \alpha_i}\right|_{\hat{\alpha}=\mathbf{1}}. \label{eq:degree}
\end{eqnarray}
The above sum considers each and every node from which an $x_i$ hyperstub originates. However, in hyperstub configuration model networks there is usually more than one type of hyperstub and this adds an additional level of detail to the excess degree. The excess degree now may incorporate two different hyperstubs into its calculations. It is now possible to describe a nodes $x_i$ degree but conditional on it being selected through ones of its $x_j$ hyperstubs. More formally we can compute the expected excess degree using conditional expectation, $E(x_j|x_i=y_i) $, which yields 
\begin{eqnarray}
\delta_{{x_j},{x_i}} = \frac{ \sum_{\hat{y}=0}^{\infty} y_j y_i p_{\hat{y}}}{ \sum_{\hat{y}=0}^{\infty} y_i p_{\hat{y}}}, \nonumber\\
\end{eqnarray}
where $\delta_{{x_j},{x_i}}$ denotes the expected $x_i$ hyperstub degree observed from a node selected proportionally to its $x_j$ hyperstub degree. The denominator is given by Eq.~\eqref{eq:Ahpgf}, and the numerator is specified by
\begin{eqnarray}
\sum_{\hat{y}^*=0}^{\infty} y_i y_j p_{\hat{y}^*} =  \left.\frac{\partial^2 \psi}{\partial \alpha_i \alpha_j}\right|_{\alpha=\mathbf{1}}. \nonumber
\end{eqnarray}

\subsection{ODEs for an example network}\label{sec:exampleodes}
The following provides ODEs for a simple example network composed of only $G_0$ and $G_\triangle$. 

When deriving ODEs by hand listing out equations for $T_i$ is a good starting point as they include many of the subgraph states, i.e. $G_0(SI)$, and can be used as the start of a check list when listing state equations  
\begin{eqnarray}
T_1 &=& \tau[G_0(SI)],\nonumber \\
T_2 &=& \tau[G_0(IS)],\nonumber\\
T_3 &=& \tau[G_\triangle(SSI) + G_\triangle(SIS) + 2G_\triangle(SII)\nonumber \\
&& +G_\triangle(SRI) + G_\triangle(SIR)],\nonumber\\
T_4 &=& \tau[G_\triangle(SSI) + G_\triangle(ISS) + 2G_\triangle(ISI)\nonumber \\
&& +G_\triangle(RSI) + G_\triangle(ISR)],\nonumber\\
T_5 &=& \tau[G_\triangle(ISS) + G_\triangle(SIS) + 2G_\triangle(IIS)\nonumber \\
&& +G_\triangle(IRS) + G_\triangle(RIS)].\nonumber
\end{eqnarray} 
It is important to node the above equations will not list every subgraph state and that for a subgraph composed of $n$ will have $3^n$ state equations. For example, the first few state equation for $G_0$ are given by
\begin{eqnarray*}
\dot{G}_0(SS) &=& -[(T \Delta)_2+ (T \Delta)_1] G_0(SS), \\
\dot{G}_0(SI) &=& -(\tau + \gamma)G_0(SI) - (T \Delta)_1 G_0(SI) +(T \Delta)_2 G_0(SS),\\
\dot{G}_0(IS) &=& -(\tau + \gamma)G_0(IS) - (T \Delta)_2 G_0(IS) +(T \Delta)_1 G_0(SS),
\end{eqnarray*}
with equations for the following being omitted 
\begin{eqnarray*}
 \{\dot{G}_0(SR),~\dot{G}_0(II),~\dot{G}_0(IR),~\dot{G}_0(RS),~\dot{G}_0(RI), \dot{G}_0(RR)\},
\end{eqnarray*}
Similarly, sample ODEs for the $G_\triangle$ subgraph, taken from a system of 27 ODEs, are: 
\begin{eqnarray*}
\dot{G}_\triangle(SSS) &=& -[(T \Delta)_5+ (T \Delta)_4 +(T \Delta)_3] G_\triangle(SSS),\\
\dot{G}_\triangle(SSI) &=& -[2\tau + \gamma + (T \Delta)_4 +(T \Delta)_3]G_\triangle(SSI) \\
&&  +(T \Delta)_5 G_\triangle(SSS),\\
\dot{G}_\triangle(SIS) &=& -[2\tau + \gamma+ (T \Delta)_5+ (T \Delta)_3]G_\triangle(SIS) \\
&& +(T \Delta)_4 G_\triangle(SSS),\\
\dot{G}_\triangle(ISS) &=& -[2\tau + \gamma + (T \Delta)_5+ (T \Delta)_4]G_\triangle(ISS) \\
&& +(T \Delta)_3 G_\triangle(SSS),
\end{eqnarray*}
with equations for the following being omitted 
\begin{eqnarray*}
 &&\{\dot{G}_0(SSR),~\dot{G}_0(SII),~\dot{G}_0(SIR),~\dot{G}_0(SRS),~\dot{G}_0(SRI),~\dot{G}_0(SRR), \\
&&~\dot{G}_0(ISI),~\dot{G}_0(ISR),~\dot{G}_0(IIS),~\dot{G}_0(III),~\dot{G}_0(IIR),~ \dot{G}_0(IRS), \\
&&~\dot{G}_0(IRI), \dot{G}_0(IRR), \dot{G}_0(RSS),~\dot{G}_0(RSI),~\dot{G}_0(RSR),~\dot{G}_0(RIS),\\
&& ~\dot{G}_0(RII), \dot{G}_0(RIR), \dot{G}_0(RRS), \dot{G}_0(RRI),~\dot{G}_0(RRR)\},
\end{eqnarray*}

Each hyperstub will have a survivor function and a corresponding ODE describing its evolution, as follows:
\begin{eqnarray*}
\dot{\theta}_1 &=& -\theta_1 \frac{T_1}{M_1}, \\
\dot{\theta}_2 &=& -\theta_2 \frac{T_2}{M_2}, \\
\dot{\theta}_3 &=& -\theta_3 \frac{T_3}{M_3}, \\
\dot{\theta}_4 &=& -\theta_4 \frac{T_4}{M_4}, \\
\dot{\theta}_5 &=& -\theta_5 \frac{T_5}{M_5}.
\end{eqnarray*}
The fraction of the population that is susceptible or infected is computed by compounding $\theta_i$ into the PGF. Symbolically, this is computed by the following
\begin{eqnarray*}
\dot{S} &=& \frac{d}{dt} \psi(\hat{\theta}), \\
\dot{I} &=& - \frac{d}{dt} \psi(\hat{\theta}) - \gamma I, \\
R &=& \gamma I,
\end{eqnarray*} 
where $\psi$ is the probability generating function that generates the hyperstub degree distribution and $\hat{\theta} = (\theta_1,\theta_2,\theta_3,\theta_4,\theta_5)$ is the probability that infection via subgraphs of types one to five has not been transmitted. The total system size for this example network is given by
\begin{eqnarray*}
3^2 + 3^3 + 5 + 2  = 43,
\end{eqnarray*}
with each term in the above corresponding to $G_0$, $G_\triangle$, survivor functions and epidemic prevalence, respectively. In general, the total number of equations is given by: 
\begin{eqnarray*}
\sum_{i=1}^M 3^{|G_i|} + |G_i|  + 2,
\end{eqnarray*}
where $G_i$ denotes a subgraph, $|G_i|$ is the number of nodes in a subgraph, and $m$ is the total number of subgraphs. 

\subsection{Equivalence to previous model for complete subgraphs}\label{sec:equiv}
The PGF formulation originally proposed by Volz et al. \citep{vmclust} is equivalent to our proposed model in the case of complete subgraphs. Consider an arbitrary complete subgraph composed of $l$ nodes and a network that is composed only of this subgraph. If positions within the subgraph are labelled distinctly, $\{x_1,x_2,\dots x_l\}$, as we have done in our approach, then the PGF of such a network is given by
\begin{eqnarray}\label{eq:hpgfl}
\psi_p(\hat{\alpha}) &=& \sum_{\hat{y}=0}^{\infty} p_{\hat{y}} \prod_{i=1}^l \alpha_i^{y_i},
\end{eqnarray}
where $\hat{y} = (y_1,y_2,\dots,y_l)$. Volz et al.'s framework treats all topologically equivalent positions as one single position. Thus, in this case, the subgraph has a single label, $x$, that corresponds to a single count, $y$, and the PGF takes the following form
\begin{eqnarray}\label{eq:hpgfv}
\psi_v(\hat{\alpha}) &=& \sum_{y=0}^{\infty} p_{y} \alpha^{y}.
\end{eqnarray}
We now show how one may obtain Eq.~\eqref{eq:hpgfv} from Eq.~\eqref{eq:hpgfl}. Since both PGFs describe the same network, the rate at which our formulation allocates position $x_i$ must be $1/l$ the rate at which Volz et al.'s formulation allocates $x$. If we replace the unique position labels of Eq.~\eqref{eq:hpgfl} with a single position marker (such as in Volz et al.'s model), the following expression is obtained
\begin{eqnarray}\label{eq:hpgfeq}
\psi_p(\hat{\alpha}) &=& \sum_{\hat{y}=0}^{\infty} p_{\hat{y}} \prod_{i=1}^l \alpha^{y/l},
\end{eqnarray}
where the following substitutions, $y_i = y/l$ and $\alpha_i = \alpha$, were made so that $\alpha^y$ is the result of the above product.
Now, every time an $x_i$ is allocated, we allocate an $x$ instead. Finally, since $p_{\hat{y}}$ is a joint distribution of $l$ identically distributed independent random variables, i.e., $\hat{y} = (y/l,y/l,\dots, y/l)$, we get: 
\begin{eqnarray*}\label{eq:collapse}
\psi_p(\hat{\alpha} =  \alpha) &=& \sum_{y=0}^{\infty} p_{y}
\alpha^{y}.
\end{eqnarray*}
It is also possible to translate between the two models elsewhere in the derivation. As an example, in our approach, infection over lines is given by $T_1$ and $T_2$, as per Eq.\eqref{eq:rofg}. By summing these values, the equivalent values used in Volz et al.'s formulation may be recovered. Following our derivation, first let $G_0(SI) \equiv G_0(IS)$ and:
\begin{eqnarray*}
T_1 + T_2 = \tau G_0(SI) + \tau G_0(IS) = 2 \tau G_0(SI).
\end{eqnarray*}
Since each $G_0$ is generated from a PGF that allocates positions at rate 1/2 that of Volz et al.'s PGF, the 2 will cancel yielding $\tau G_0(SI)$. However, it is only necessary to show equivalence between the two PGFs since all other variables follow from this.

\subsection{State transition matrix}\label{sec:statetm}

The state transition matrix for $G_0$ (lines) is given by:
\newcommand\scalemath[2]{\scalebox{#1}{\mbox{\ensuremath{\displaystyle #2}}}}
\begin{center}
$\mathbf{Z}$ = \scalemath{0.88}{\bordermatrix{~    & (SS) & (SI) & (SR) & (IS) & (II) & (IR) & (RS) & (RI) & (RR)\cr
                  (SS) &   0  &  (T \Delta)_2   &     0& (T \Delta)_1 &     0&     0&     0&     0&  0  \cr
                  (SI) &   0  &  0   &   \gamma &     0&     \tau + (T \Delta)_1 &     0&     0&     0&  0  \cr
                  (SR) &   0  &  0   &     0&     0&     0&     (T\Delta)_1 &     0&     0&  0  \cr
                  (IS) &   0  &  0   &     0&     0&  \tau+(T\Delta)_2   &     0& \gamma & 0 &  0  \cr
                  (II) &   0  &  0   &     0&     0&     0& \gamma &     0& \gamma &  0  \cr
                  (IR) &   0  &  0   &     0&     0&     0&     0&     0 &     0&  \gamma  \cr
                  (RS) &   0  &  0   &     0&     0&     0&     0&     0&     (T\Delta)_2&  0  \cr
                  (RI) &   0  &  0   &     0&     0&     0&     0&     0&     0&  \gamma  \cr
                  (RR) &   0  &  0   &     0&     0&     0&     0&     0&     0&  0  \cr}}.
 \end{center}

\newpage 
 
\subsection{Algorithm 1 - Hyperstub CM algorithm}\label{sec:alg1}
\def\HiLiy{\leavevmode\rlap{\hbox to \hsize{\color{yellow!30}\leaders\hrule height .8\baselineskip depth .5ex\hfill}}}
\def\HiLig{\leavevmode\rlap{\hbox to \hsize{\color{green!20}\leaders\hrule height .8\baselineskip depth .5ex\hfill}}}
\def\HiLio{\leavevmode\rlap{\hbox to \hsize{\color{orange!50}\leaders\hrule height .8\baselineskip depth .5ex\hfill}}}
\def\HiLib{\leavevmode\rlap{\hbox to \hsize{\color{blue!20}\leaders\hrule height .8\baselineskip depth .5ex\hfill}}}
\def\HiLir{\leavevmode\rlap{\hbox to \hsize{\color{red!20}\leaders\hrule height .8\baselineskip depth .5ex\hfill}}}
\begin{algorithm}
\setlength{\leftskip}{10pt}
\SetKwData{Left}{left}\SetKwData{This}{this}\SetKwData{Up}{up}
\SetKwInOut{Input}{input}\SetKwInOut{Output}{output}
\HiLiy \Input{$N$, $K$,} 
\HiLiy \Output{$A$.}
\BlankLine
\HiLib \KwSty{Variables / initialisation} \\
$N$: the number of nodes, \\
\HiLig \CommentSty{\% Each row of $K$ corresponds to single node's hyperstub sequence.} \\
$K$: the hyperstub degree sequence, a non-square matrix $K \in \mathbb{N}_0^{N \times H}$,\\
$H$: the number of hyperstub types, \\
$A$: the adjacency matrix of the network, $A \in \{0,1\}^{N \times N}$, \\
$M$: the number of subgraphs, \\
$H_i$: the degree of a specific hyperstub, \\
$h_i$: a dynamic list of nodes that are incident to $H_i$ hyperstubs, \\
$g_i$: the adjacency matrix of a subgraph, $g \in \{0,1\}^{n_i \times n_i}$, \\
$n_i$: the number of nodes in $g_i$. \\
\BlankLine
\HiLio \KwSty{Procedure} \\
\HiLig \CommentSty{\% The following creates dynamic lists the, `hyperstub bins'.}  \\ 
\For{\DataSty{every node $i$}}{
	\For{\DataSty{each $H_j$}}{
		\DataSty{append $K_{i,j}$ multiples of $node(i)$ to the hyperstub bin($h_j$)} \\
	}}
\For{\DataSty{For each subgraph $g_i$}}{
	\For{\DataSty{For each hyperstub of $g_i$}}{
\HiLig \CommentSty{\% Select unfiformly at random and without} \\
\HiLig \CommentSty{\% replacement a node incident to each desired hyperstub.}  \\ 
		\DataSty{$n_1$ = rand-sample($h_{i_1}$)} \label{1st:node1} \\
		\DataSty{$n_2$ = rand-sample($h_{i_2}$)} \label{1st:node2} \\
		\DataSty{~~~~ \vdots} \\
		\DataSty{$n_{g_i}$ = rand-sample($h_{i_2}$)} \label{1st:noden} \\
	}
	\HiLig \CommentSty{\% The following compares pairs of the selected nodes} \\
	\HiLig \CommentSty{\% to determine their connectivity in $A$. } \\
	\For{\DataSty{$k = (1,2,\dots,n_i)$}}{
			\For{\DataSty{$l = (1,2,\dots,n_i)$}}{
		\DataSty{\If{$g(n_k,n_l)==1$}{
	\DataSty{$A(n_k,n_l)=1$ } \\\label{1st:over}
	}}}}
}\caption{The hyperstub configuration model. In this implementation, multiple-edges are over written (line~\ref{1st:over}) but self-edges are permitted. To prevent this, if nodes already share an edge or a node has been selected twice (self-edge) lines~\ref{1st:node1}-\ref{1st:noden} are repeated until a valid selection is made. This reselection step has been omitted below for readability.}\label{algo:HCMA}
\end{algorithm}

\newpage

\subsection{Algorithm 2 - Transition matrix algorithm}\label{sec:alg2}

\begin{algorithm}[ht!]
\SetKwData{Left}{left}\SetKwData{This}{this}\SetKwData{Up}{up}
\SetKwFunction{Union}{Union}\SetKwFunction{FindCompress}{FindCompress}
\SetKwInOut{Input}{input}\SetKwInOut{Output}{output}
\HiLiy\Input{$g$,}
\HiLiy\Output{$\mathbf{Z}$.}
\BlankLine
\HiLib\KwSty{Variables / initialisation} \\
$g$: the adjacency matrix of a subgraph $G$,\\
\HiLig \% $\mathbf{Z} \in \mathbb{R}^{3^n \times 3^n}$. \\
$\mathbf{Z}$: matrix corresponding rate of transition between states of $G$, \\
$n$: node count of $G$, \\
\HiLig \% $\vv{G}$ contains $3^n$ elements. \\
$\vv{G}$: the vector of states of $G$, \\
$\tau$: per link infection rate, \\
$\gamma$: recovery rate, \\
$ T\Delta$: the expected force of infection a node within $G$ experiences from outside $G$. 
\BlankLine
\HiLio\KwSty{Procedure} \\
\For{every state $\vv{G}_i$}{
\For{every state $\vv{G}_j$}{
\HiLig \% Compare each and every possible state transition of $G$:\\ 
\Switch{ $\vv{G}_i \rightarrow  \vv{G}_j$ }{\label{2nd:switch}
\Case{A single infection occurs}{
\eIf{the new $I$ is connected to another $I$ within $G$}{
\HiLig \% Check the connectivity of the new $I$ using $g$. \\
 $Z_{i,j} = \tau + T\Delta$}
{\HiLig \% the infection was from only an external source. \\
$Z_{i,j} = T\Delta$}}
\Case{A single recovery occurs}{
$Z_{i,j} = \gamma$}
\Case{otherwise}{
$Z_{i,j} = 0$}
}}}\caption{Generating the state transition matrix. The comparison in line~\ref{2nd:switch} needs to check: (1) that only a single node has changed state and (2) only state changes $S \rightarrow I$ and $I \rightarrow R$ are valid.}\label{algo:TMA}
\end{algorithm}

 \newpage

\subsection{Null case for Fig.~\ref{fig:allplotsc}}\label{sec:null}

\begin{figure}[ht!]
\begin{center}
 \includegraphics[scale=0.6]{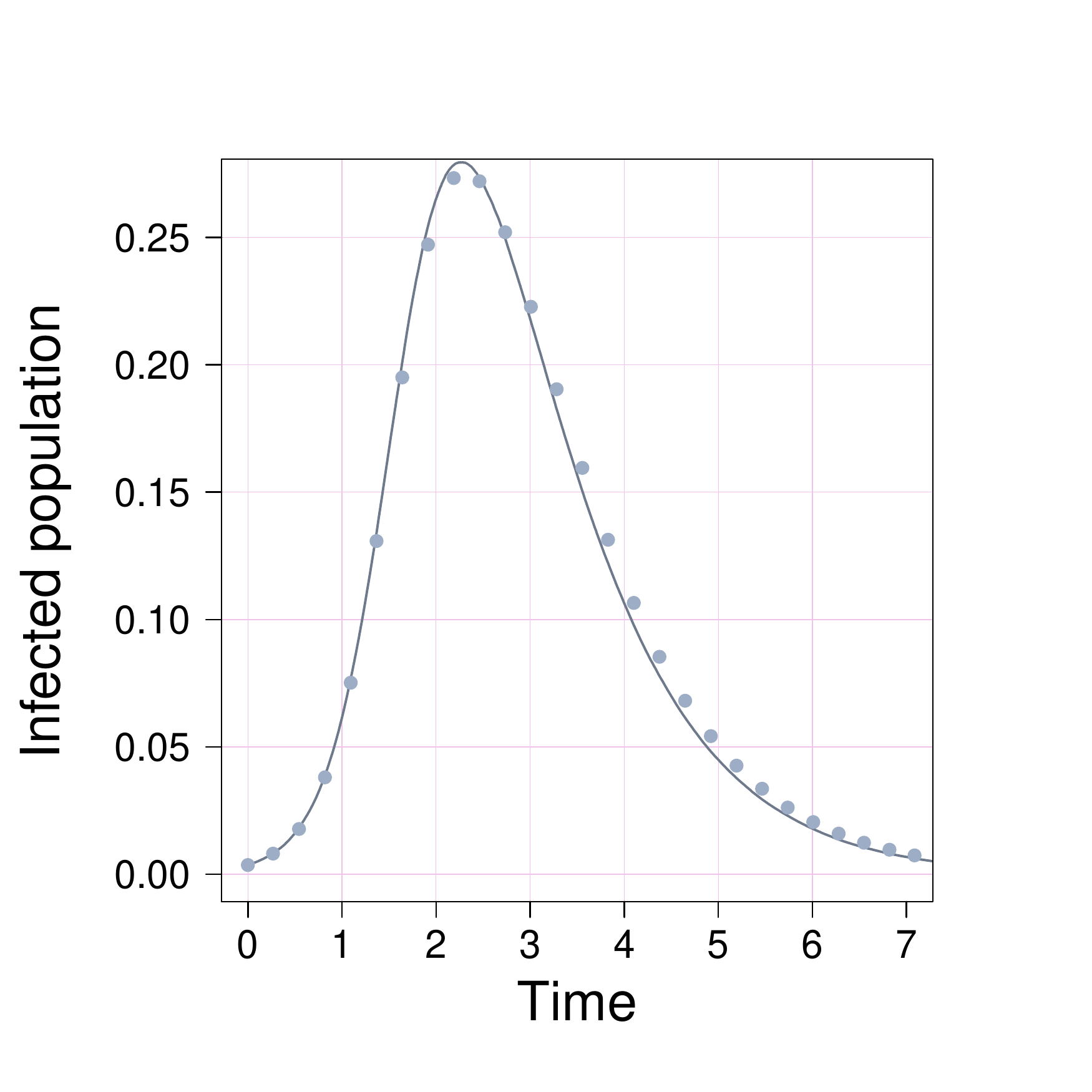}
 \end{center}
 \caption{The effect of higher moments. The solid and discrete plots correspond to the null networks $G_0 \sim 2Pois(2)$ and $G_0 \sim Pois(3) + 5Pois(1/5)$ respectively, i.e. the null cases for the triangle and hexagon networks. Both plots have equal first and second moments and clustering equal to that of a random network. The difference observed is a result of non-equal higher moments and is not enough to explain the difference observed in Fig.~\ref{fig:allplotsc}.}
 \label{fig:highmom}
 \end{figure}

\end{document}